\documentclass{amsart}
\newcommand{\dbar}{\bar\partial}
\newcommand{\Cn}{\mathbb{C}^n}
\newcommand{\C}{\mathbb{C}}
\newcommand{\CN}{\mathbb{C}^N}
\newcommand{\CNp}{\mathbb{C}^{N^\prime}}
\newcommand{\Rd}{\mathbb{R}^d}
\newcommand{\Rn}{\mathbb{R}^n}
\newcommand{\RN}{\mathbb{R}^N}
\newcommand{\R}{\mathbb{R}}
\newcommand{\N}{\mathbb{N}}
\newcommand{\da}{\delta}
\newcommand{\Om}{\Omega}
\newcommand{\spaceb}{\mathfrak{B} (\Om_+)}
\newcommand{\spacebm}{\mathfrak{B} (\Om_-)}
\newcommand{\spacea}{\mathfrak{A} (\Om_+)}
\newcommand{\spaceainf}{\mathfrak{A}_\infty (\Om_+)}
\newcommand{\spaceam}{\mathfrak{A} (\Om_-)}
\newcommand{\spaceaminf}{\mathfrak{A}_\infty (\Om_-)}
\newcommand{\almhol}{\mathfrak{AH} (U\times V)}
\newcommand{\dop}[1]{\frac{\partial}{\partial #1}}
\newcommand{\vardop}[2]{\frac{\partial #1}{\partial #2}}
\newcommand{\br}[1]{\langle#1 \rangle}
\newcommand{\eps}{\epsilon}
\DeclareMathOperator{\supp}{supp}
\DeclareMathOperator{\spanc}{span}
\DeclareMathOperator{\imag}{Im}
\newtheorem{thm}{Theorem}
\newtheorem{lem}[thm]{Lemma}

\newtheorem{defin}[thm]{Definition}
\newtheorem{cor}[thm]{Corollary}
\newtheorem{cl}{Claim}
\begin{document}
\title[$C^\infty$-regularity for  nondegenerate CR-mappings]{A $C^\infty$-regularity theorem for\\ nondegenerate {C}{R} mappings}%
\author{Bernhard Lamel }%
\address{Kungliga Tekniska H\"ogskolan, Stockholm}%
\email{lamelb@member.ams.org}%
\thanks{The author would like to thank the Erwin Schr\"odinger
Institut in Vienna for its hospitality.}

\begin{abstract} We prove the following regularity result: If
$M\subset \CN$, $M^\prime\subset\CNp$ are smooth generic
submanifolds and $M$ is minimal, then every $C^k$-CR-map from $M$
into $M^\prime$ which is $k$-nondegenerate is smooth. As an
application, every CR diffeomorphism of $k$-nondegenerate minimal
submanifolds  in $\CN$ of class $C^{k}$ is smooth.
\end{abstract}
\maketitle
\section{Introduction and statement of results}
We first briefly describe the setting for the results which we
want to discuss. Let $M\subset\CN$, $M^\prime\subset\CNp$ be
generic, real submanifolds of $\CN$ and $\CNp$, respectively. We
shall denote by $d$ the real codimension of $M$ and by $d^\prime$
the real codimension of $M^\prime$, and write $n=N-d$, $n^\prime =
N^\prime - d^\prime$. Recall that $M$ is generic if there  is a
smooth defining function $\rho = (\rho_1,\dots,\rho_d)$ for $M$
such that the vectors $\rho_{1,Z}(p),\dots, \rho_{d,Z}(p)$ are
linearly independent for $p\in M$. Here for any smooth function
$\phi$ we let $\phi_Z =
(\vardop{\phi}{Z_1},\dots,\vardop{\phi}{Z_N})$ be its complex
gradient.

We also fix points $p_0\in M$ and $p_0^\prime\in M^\prime$ (which
we will assume to be equal to $0$ for most of this paper). A
$C^k$-mapping $H$ from $M$ into $M^\prime$ is said to be CR if
its differential $dH$ satisfies $dH(T^c_p M)\subset T^c_{H(p)}
M^\prime$ for $p\in M$, where $T^c_p M$ denotes the complex
tangent space to $M$ at $p$, that is, the largest subspace of the
real tangent space $T_p M$ invariant under the complex structure
operator $J$ in $\CN$. Equivalently, if $H=(H_1,\dots,
H_{N^\prime})$ for any system of holomorphic coordinates in
$\CNp$, each $H_j$ is a CR-function on $M$. (For further
reference on these definitions, the reader is referred to the
book of Baouendi, Ebenfelt and Rothschild \cite{BERbook}).

The following definition is from \cite{L1}. We shall give it in a
slightly modified form.

\begin{defin}\label{D:nondeg} Let $M$, $M^\prime$ be as above. Let $\rho^\prime =
(\rho_1^\prime,\dots,\rho_{d^\prime}^\prime)$ be a defining
function for $M^\prime$ near $H(p_0)$, and choose a basis
$L_1,\dots, L_n$ of CR-vector fields tangent to $M$ near $p_0$.
We shall write $L^\alpha = L_1^{\alpha_1}\cdots L_n^{\alpha_n}$
for any multiindex $\alpha$. Let $H:M\to M^\prime$ be a CR-map of
class $C^m$. For $0\leq k \leq m$, define the increasing sequence
of subspaces $E_k(p_0)\subset\CNp$ by
\begin{equation} \label{E:theeks}
  E_k (p_0) = \spanc_\C  \{ L^\alpha \rho_{l,Z^\prime}^\prime
  (H(Z),\overline{H(Z)})|_{Z=
  p_0} \colon
  0\leq |\alpha|\leq k, 1\leq l \leq d^\prime\}.
\end{equation}
We say that $H$ is $k_0$-nondegenerate at $p_0$ (with $0\leq k_0
\leq m$) if $E_{k_0 -1}(p_0)\neq  E_{k_0} (p_0) = \CNp$.
\end{defin}
The invariance of this definition under the choices of the
defining function, the basis of CR vector fields and the choices
of holomorphic coordinates in $\CN$ and $\CNp$ is easy to show;
the reader can find proofs for this in \cite{L1} or \cite{L2}.

Recall that if $\Gamma\subset\Rd$ is an open convex cone, $p_0
\in M$, and $U\subset \CN$ is an open neighbourhood of $p_0$,
then a wedge $W$ with edge $M$ centered at $p_0$ is defined to be
a set of the form $W=\{ Z\in U \colon \rho(Z,\bar Z)
\in\Gamma\}$, where $\rho$ is a local defining function for $M$.
We can now state our main theorem.
\begin{thm}\label{T:main}
Let $M\subset\CN$, $M^\prime\subset\CNp$ be smooth generic
submanifolds of $\CN$ and $\CNp$, respectively, $p_0\in M$ and
$p_0^\prime\in M^\prime$, $H:M\to M^\prime$ a $C^{k_0}$-CR-map
which is $k_0$-nondegenerate at $p_0$ and extends continuously to
a holomorphic map in a wedge $W$ with edge $M$. Then $H$ is smooth
in some neighbourhood of $p_0$.
\end{thm}
This theorem is the smooth version of the main result in
\cite{L1}. Let us recall that $M$ is said to be minimal at $p_0$
if there does not exist any CR-submanifold through $p_0$ strictly
contained in $M$ with the same CR dimension as $M$. By a theorem
of Tumanov, if $M$ is minimal, every continuous CR-function $f$
on $M$ near $p$ extends continuously to a holomorphic function
into a wedge $W$ with edge $M$. Hence we have the following
corollary.

\begin{cor} Let $M\subset\CN$, $M^\prime\subset\CNp$ be smooth generic
submanifolds of $\CN$ and $\CNp$, respectively, $p_0\in M$ and
$p_0^\prime\in M^\prime$, $M$ minimal at $p_0$, $H:M\to M^\prime$
a $C^{k_0}$ map which is $k_0$-nondegenerate at $p_0$. Then $H$
is smooth in some neighbourhood of $p_0$.
\end{cor}

Note that by a regularity theorem of Rosay (\cite{Ro1}, see also
\cite{BERbook}), if the boundary value of a holomorphic function
in a wedge $W$ with edge $M$ is $C^k$ on $M$, then the extension
is also of class $C^k$ up to the edge. Hence, for the proof of
Theorem~\ref{T:main} we will assume that $H$ extends in a
$C^{k_0}$-fashion to a wedge $W$ centered at $p_0$.

We would like to mention one particular instance of this theorem.
If $M$ is a manifold whose identity map is $k_0$-nondegenerate in
the sense of Definition~\ref{D:nondeg}, then we say that $M$ is
{\em $k_0$-nondegenerate}. This notion has been introduced for
hypersurfaces by Baouendi, Huang and Rothschild in \cite{BHR1};
for a thorough introduction to this nondegeneracy condition for
submanifolds and its connection with holomorphic nondegeneracy in
the sense of Stanton (\cite{S2}), see \cite{BERbook}, or the
paper of Ebenfelt \cite{E3}. In particular, every
CR-diffeomorphism of class $C^{k_0}$ of a $k_0$-nondegenerate
submanifold is $k_0$-nondegenerate in the sense of
Definition~\ref{D:nondeg}. Theorem~\ref{T:main} implies the
following regularity result for $k_0$-nondegenerate smooth
submanifolds.

\begin{cor}\label{C:nondeg} Assume that $M\subset\CN$ and $M^\prime\subset\CNp$ are
$k_0$-nondegenerate smooth submanifolds of real codimension $d$,
$M$ minimal at $p_0$, and $H:M\to M^\prime$ is a
CR-diffeomorphism of class $C^{k_0}$. Then $H$ is smooth.
\end{cor}

If $d=1$, we can drop the assumption of minimality, since in the
hypersurface case, $k_0$-nondegeneracy implies minimality. In the
case where $N=N^\prime=2$ and $d=1$, Corollary~\ref{C:nondeg} is
basically contained in the thesis of Roberts \cite{Rob}. The
Levi-nondegenerate hypersurface case is well understood; the
connection with the results proved in this paper is that
Levi-nondegeneracy of hypersurfaces is equivalent to
$1$-nondegeneracy. In fact, for Levi-nondegenerate hypersurfaces,
Corollary~\ref{C:nondeg} is due to Nirenberg, Webster and Yang
\cite{NWY}, and of course we should not forget to mention
Fefferman's mapping theorem \cite{Fe1} (however, we shall not
deal with the $C^1$-extension here). A proof for strictly
pseudoconvex hypersurfaces of finite smoothness was given by
Pinchuk and Khasanov \cite{PK}. More recently, Tumanov
\cite{Tu94} has proved the corresponding theorem for
Levi-nondegenerate targets of higher codimension. For results for
pseudoconvex targets, we want to refer the reader to the
historical discussion in the paper by Coupet and Sukhov
\cite{CS97} and the newer results for convex hypersurfaces by
Coupet, Gaussier and Sukhov \cite{CS99}.

The paper is organized as follows. In
section~\ref{S:boundaryvalues}, \ref{S:edgeofthewedge}, and
\ref{S:ift} we present the technical foundations for the proof.
Although these results are well known, they are not easy to find
in the literature; so, in order to make this paper as self
contained as possible, we have decided to include the proofs.
Theorem~\ref{T:main} is then proved in section~\ref{S:proof}.

\section{Boundary values of functions of slow
growth}\label{S:boundaryvalues} In this section, we will develop
an integral representation for a $\dbar$-bounded function of slow
growth (in a wedge with straight edge). Let us first fix
notation. Let $U\subset\Cn$, $V\subset\Rd$ be open subsets, and
let $\da = (\da_1,\dots,\da_d)\in\Rd$ with $0<\da_j$ for $0\leq
j\leq d$. We set $\Om_+ = \{(z,s,t)\in U\times V\times \Rd \colon
0<t<\da \}$, $\Om_- =\{(z,s,t)\in U\times V\times \Rd \colon
0>t>-\da \}$ and $\Om_0 = U\times V\times \{0 \}$, and we will
write $z=(x,y)$ for the underlying real variables. Throughout the
paper, $dm$ will denote Lebesgue measure. Let $\spaceb$ be the
space of all functions $h\in C^1(\Om_+)$ that extend smoothly to
the set $E = \{ (z,s,t) \in \bar\Om_+ \colon t\neq 0\}$ which
have the following property: For each compact set $K\subset
U\times V$, there exist positive constants $C_1$, $\mu$ and $C_2$
(depending on $K$ and $h$) such that
\begin{equation}\label{E:slowgrowth}
\sup_{(z,s)\in K,0<t<\da  }|t|^\mu |h(x,y,s,t)| \leq C_1
\end{equation}
and
\begin{equation}\label{E:dbarbounded}
\sup_{(z,s)\in K,0<t<\da  }|\dbar_j h(x,y,s,t)| \leq C_2, \quad
1\leq j\leq d.
\end{equation}
Here we write $\dbar_j= \frac{1}{2} \left(\dop{s_j} + i \dop{t_j}
\right)$. We have the following (probably well known) result,
which we state for $\spaceb$; however, we  define $\spacebm$ in a
similar manner, and all the results stated in this section hold
equally well for $\spacebm$.
\begin{thm}\label{T:boundaryvalue}  Let $h\in \spaceb$. Then the limit
\begin{equation}\label{E:boundaryvaluedefined}
  \br{ b_+ h , \phi } = \lim_{\eps = (\eps_1,\dots ,\eps_d) \to
  0} \int_{U\times V} h(x,y,s,\eps) \phi (x,y,s) \, dm
\end{equation}
exists for each $\phi \in C^\infty_c (\Om_0 )$ and defines a
distribution $b_+ h$ called the boundary value of $h$.
Furthermore, for each compact set $K$  there exists an integer
$v_0$ such that for $v\geq v_0$, for each $j=1,\dots, d$, $0\leq
\delta^\prime \leq \delta_j$ we have the following integral
representation for $\phi\in C^\infty_c(U\times V)$ with $\supp
\phi \subset K$:
\begin{multline}\label{E:intrep}
\br{b_+h, \phi} = \int_{U\times V}
h(x,y,s,0,\dots,\da^\prime,\dots,0) S_v \phi (x,y,s,0,\dots,
\delta^\prime,\dots, 0) \,
dm \\
+ 2i \int_{U\times V} \int_0^{\da^\prime}  \dbar_1 h
(x,y,s,0,\dots,t_j,\dots,0) S_v\phi (x,y,s,0,\dots,t_j,\dots,0) \,  dt_j dm\\
+ 2i  \int_{U\times V}  \int_0^{\da^\prime}  h (x,y,s,0,\dots,
t_j,\dots,0 ) D_{s_j}^{v+1} \phi (x,y,s) t_j^v \, dt_j dm.
\end{multline}
where
\begin{equation}\label{E:definsv}  S_v \phi (x,y,s,t) =
\sum_{|\alpha| \leq v} {\frac{1}{\alpha!}} D^\alpha_s \phi(x,y,s)
t^\alpha.
\end{equation}
\end{thm}
\begin{proof}
Let $S_v \phi$ be defined by \eqref{E:definsv}. We are going to
prove the formula under the assumption that $j=1$. Fix $(x,y)$,
$s_2,\dots s_d$ and $0<\da^\prime<\da_1$, and assume
$0<\epsilon_1 <\da_1 - \da^\prime$. First we are going to assume
that $K=\supp \phi$ is contained in a product of the form
$U_1\times[a,b]\times [a_2,b_2]\times \dots \times [a_d, b_d]$
contained in a relatively compact open subset $W\subset U\times
V$. In this case, define
\[u(s_1,t_1) = h(x,y,s_1,s_2,\dots, s_d, \eps_1 +
t_1,\eps_2,\dots ,\eps_d) S_v \phi (z,s,t_1,0,\dots, 0).\]
Clearly, $u$ is $C^1$ on the square $\omega =
[a,b]\times[0,\da^\prime]$ and $u(s_1,t_1) = 0$ if $s_1\geq b$
or  $s_1\leq a$. By Stokes formula,
\[ \int_{\partial \omega} u(s_1, t_1) \, dw = 2i \int_\omega \dbar u
(s_1,t_1) \, dm ,\] where we have set $w = s_1 + i t_1$ and
$\dbar = \dbar_1$. This formula translates into
\begin{multline} \int_a^b h(x,y,s,\eps) \phi(x,y,s) \, ds_1 =\\
\int_a^b h(x,y,s, \eps_1+\da^\prime, \eps_2,\dots, \eps_d) S_v
\phi (x,y,s,\delta^\prime,0,\dots, 0) \,
ds_1 \\
+ 2i \int_0^{\da^\prime} \int_a^b \dbar_1 h (x,y,s, \eps_1+ t_1,
\eps_2,\dots, \eps_d ) S_v\phi (x,y,s,t_1,0,\dots,0) \, ds_1 dt_1 \\
+ 2i \int_0^{\da^\prime} \int_a^b \dbar_1 h (x,y,s, \eps_1+ t_1,
\eps_2,\dots, \eps_d ) D_{s_1}^{v+1} \phi (x,y,s) t_1^v \, ds_1
dt_1.
\end{multline}
We integrate this formula with respect to  $(x,y,s_2,\dots, s_d)$
to obtain
\begin{multline}
\int_{W} h(x,y,s,\epsilon) \phi (x,y,s) \, dm =\\
\int_{W} h(x,y,s, \eps_1+\da^\prime, \eps_2,\dots, \eps_d)S_v \phi
(x,y,s,\delta^\prime,0,\dots, 0) \,
dm \\
+ 2i \int_{W} \int_0^{\da^\prime}  \dbar_1 h (x,y,s, \eps_1+ t_1,
\eps_2,\dots, \eps_d ) S_v\phi (x,y,s,t_1,0,\dots,0) \,  dt_1 dm\\
+ 2i  \int_{W}  \int_0^{\da^\prime}  h (x,y,s, \eps_1+ t_1,
\eps_2,\dots, \eps_d ) D_{s_1}^{v+1} \phi (x,y,s) t_1^v \, dt_1
dm.
\end{multline}
For each of these integrals, we can use the bounded convergence
theorem to take the limit as $\eps \to 0$, provided that we
choose $v\geq \mu_K$, where $\mu_K$ denotes the least integer
$\mu$ for which \eqref{E:slowgrowth} holds on $K$  and to obtain
an estimate of the form $|\br{b_+ h, \phi}| \leq C\| \phi
\|_{v+1}$ (where $\| \phi \|_{k} = \max_{x\in U\times V,
|\alpha|\leq k} |\phi^\alpha (x)|)$.

Now we pass to the case of general $K$ by covering with finitely
many sets of the form considered above and using a partition of
unity. The details are easy and left to the reader.
\end{proof}

Consider now the class $\spacea$ of functions $h$ which are smooth
on $E$ with the property that for all $\alpha,\beta$ we have that
$D^{\alpha}_{x,y} D^\beta_s h\in \spaceb$. If $h\in\spacea$, for
$K\subset U\times V$ we let $\mu_l (h,K)$ the smallest integer
$\mu$ such that
\begin{equation}\label{E:slowgrowtha}
\sup_{(z,s)\in K,0<t<\da }|t|^\mu |D^{\alpha}_{x,y} D^\beta_s
h(x,y,s,t)| \leq C_1, \quad |\alpha|+|\beta|\leq l
\end{equation}
for some constant $C_1$. Let us also introduce the space
$\spaceainf $ of functions in $\spacea$ with the additional
property that for any compact set $K\subset U\times V$,  for any
multiindeces $\alpha$ and $\beta$, and for any nonnegative
integer $k$ there exists a constant $C$ such that
\begin{equation}\label{E:dbartozero}
\sup_{(z,s)\in K,0<t<\da  }|D^{\alpha}_{x,y} D^\beta_s \dbar_j
h(x,y,s,t)| \leq C |t|^k, \quad 1\leq j\leq d.
\end{equation}
Of course, we define the spaces $\spaceam$ and $\spaceaminf$
analogously, and the results stated below for $\spacea$ and
$\spaceainf$ also hold for $\spaceam$ and $\spaceaminf$. This can
be seen most easily by noting the following useful fact: If
$h(x,y,s,t) \in \spacea$ (or $\spaceainf$, respectively),
$\overline{h(x,y,s,-t)}\in \spaceam$ (or $\spaceaminf$,
respectively).

We will also need the space of functions which are {\em almost
holomorphic} on $U\times V$. This is the space
\begin{equation}\label{E:almhol}
\almhol = \{ a\in C^\infty (U\times V\times \Rd) \colon
D^\alpha_{x,y} D^\beta_s D^\gamma_t \dbar_j a(x,y,s,0) = 0, 1\leq
j\leq d \}.
\end{equation}
\begin{lem} \label{L:multlemma} Let $h\in\spacea$, $a\in\almhol$,
and set $a_0(x,y,s) = a(x,y,s,0)$. Then $ah\in\spacea$, and $b_+
ah = a_0 b_+h$ in the sense of distributions. Furthermore, if
$h\in \spaceainf$, so is $ah$.
\end{lem}
\begin{proof} By the Leibniz rule, $D^{\alpha}_{x,y} D^\beta_s ah$
is a sum of products of derivatives of $a$ and $h$. It is clear
that such a sum fulfills \eqref{E:slowgrowth}. To see that it
also fulfills \eqref{E:dbarbounded}, note that by
\eqref{E:almhol} every derivative of $\dbar_j a$ vanishes to
infinite order on $t=0$.

To see that $b_+ ah = a_0 b_+ h$ we use Taylor development to
write $a(x,y,s,t)  = \sum_{|\beta|\leq k} \frac{1}{\beta !}
D^\beta_s a(x,y,s,0) (it)^\beta + O(|t|^{k+1})$ (uniformly on
compact subsets of $U\times V$). Now choose $k\geq \mu_0 (h,K)$
and substitute into \eqref{E:boundaryvaluedefined} for $\phi$
with $\supp \phi \subset K$. The claim follows now by taking the
limit and using Theorem~\ref{T:boundaryvalue}.
\end{proof}
Basically the same proof shows the following Lemma.
\begin{lem}\label{L:difflemma} Assume that $X$ is a vector field
on $U\times V \times \Rd$ which is tangent to all subspaces of
the form $t = c$, where $c\in\Rd$ is a constant vector, and such
that all the coefficients of $X$ are in $\almhol$. Set $X_0 =
X|_{t=0}$. If $h\in\spacea$, then $Xh\in\spacea$, and $b_+ Xh =
X_0 b_+ h$ in the sense of distributions. Furthermore, if
$h\in\spacea$, so is $Xh$.
\end{lem}

\section{An almost holomorphic edge-of-the-wedge
theorem}\label{S:edgeofthewedge} The main result of this section
is the following theorem. Our presentation follows closely
\cite{Rob}, but we also want to refer the reader to \cite{Ros2}.
We keep the notation from the proceeding section and since we
shall use the Fourier transform we also introduce the following
new variables: $\xi \in \Rn$, $\tau\in\Rn$, $\sigma\in\Rd$. For a
distribution $\phi$ on $U\times V$ we will write $\hat \phi (\xi,
\tau, \sigma) = \br{ \phi, \exp(-i(x\xi + y\tau + s\sigma))}$ for
its Fourier transform.

\begin{thm}\label{T:edgeofthewedge} Assume that $h_+ \in \spacea$,
$h_-\in\spaceam$, and
that $b_+ h_+=b_- h_- = h$. Then $h$ is smooth.
\end{thm}

The proof follows from the next Lemma.

\begin{lem} \label{L:estabove} Let $h \in \spacea$, and
$\phi\in C^\infty_c (U\times V)$.
 Then for every $k\in\N$ there
exists a constant $C_k$ such that if $\zeta=
(\xi,\tau,\sigma)\in\Rn\times\Rn\times\Rd$ with $\sigma_j\leq 0$
for some $j$, $1\leq j \leq d$, then
\begin{equation}\label{E:estfork}
| \widehat{\phi b_+ h}
 (\zeta) | \leq \frac{C_k}{(1+|\zeta |^2)^k}.
\end{equation}
Here, $C_k$ depends on $k$, $\phi$, and $h$. The same result holds
with $\spacea$ replaced by $\spaceam$ if $\sigma_j \geq 0 $ for
some $j$.
\end{lem}

\begin{proof} For the moment, fix $\zeta$; for simplicity,
assume that $j=1$, so that $\sigma_1 \leq 0$. We shall write
$a(x,y,s,t) = \exp( -i(x\xi + y\tau + s\sigma) + t\sigma)$. Then
$a\in\almhol$---in fact, $\partial_j a = 0$, $1\leq j\leq d$. We
let $\Delta$ be the real Laplacian in the $2n+d$ variables
$(x,y,s)$, that is,
\begin{equation}\label{E:laplace}
\Delta = \sum_{j=1}^n \vardop{^2}{x_j^2} + \sum_{j=1}^n
\vardop{^2}{y_j^2} + \sum_{j=1}^d \vardop{^2}{s_j^2}.
\end{equation}
We then have that $(1+\Delta)^k a(x,y,s,t) = (1+|\zeta|^2)^k
a(x,y,s,t)$. Recall that we write $a_0(x,y,s) = a(x,y,s,0)$. By
Lemma~\ref{L:multlemma}, we see that $\widehat{\phi b_ +h} (\zeta)
= \br{\phi b_+ h, a_0} = \br{ b_+ h , \phi a_0} = \br{ a_0 b_+ h,
\phi} = \br{ b_+ ah , \phi}$. We apply the integral formula
\eqref{E:intrep} from Theorem~\ref{T:boundaryvalue} for $j = 1$,
and some $\da^\prime$, which implies that
\begin{multline}\label{E:intrepapp}
\br{b_+ a h, \phi} = \\
\int_{U\times V} h(x,y,s,\da^\prime,0) e^{-i(x\xi + y\tau +
s\sigma)}e^{\da^\prime \sigma_1} S_v \phi (x,y,s, \delta^\prime,
0) \,
dm \\
+ 2i \int_{U\times V} \int_0^{\da^\prime}  (\dbar_1 h
(x,y,s,t_1,0)) e^{-i(x\xi + y\tau +
s\sigma)}e^{t_1 \sigma_1} S_v\phi (x,y,s,t_1,0) \,  dt_1 dm\\
+ 2i  \int_{U\times V}  \int_0^{\da^\prime}  h (x,y,s, t_1,0 )
e^{-i(x\xi + y\tau + s\sigma)}e^{t_1 \sigma_1} D_{s_1}^{v+1} \phi
(x,y,s) t_1^v \, dt_1 dm\\
= I_1 + I_2 + I_3.
\end{multline}
We now replace $e^{-i(x\xi + y\tau + s\sigma)}$ by
$\frac{1}{(1+|\zeta|^2)^k} (1+ \Delta)^k e^{-i(x\xi + y\tau +
s\sigma)}$ in all three integrals above. Then we integrate by
parts and estimate, where we choose $v\geq \mu_{2k} (h,K)$ (see
\eqref{E:slowgrowtha} for the definition of this number) with
$K=\supp \phi$. Since all the estimates are easy, we do not write
them out; the reader can easily check them.
\end{proof}
\begin{proof}[Proof of Theorem~\ref{T:edgeofthewedge}] Let $p\in U\times
V$. Choose a function $\phi\in C^\infty_c (U\times V)$ which is
equal to $1$ in some open neighbourhood of $p$. By
Lemma~\ref{L:estabove}, since $h_+\in\spacea$ and
$h_-\in\spaceam$, we have that
\begin{equation}\label{E:estfork2}
| \widehat{\phi h}
 (\zeta) | \leq \frac{C_k}{(1+|\zeta |^2)^k}.
\end{equation}
for all $\zeta\in \R^{2n+d}$. Hence, $\phi h$ is smooth (see for
Example \cite{Foll}), and so $h$ is smooth in some neighbourhood
of $p$, since $\phi \equiv 1$ there. Since $p$ was arbitrary, the
claim follows.
\end{proof}
\section{A version of the Implicit Function Theorem}\label{S:ift}
We will need the following, ``almost holomorphic'', implicit
function theorem.
\begin{thm} \label{T:ift} Let $U\subset\CN$ be open, $0\in U$,
$A\in \C^p$, $F:U\times \C^p \to\CN$ be smooth in the first $N$
variables and polynomial in the last $p$ variables, and assume
that $F(0,A) = 0$ and $F_Z (0,A)$ is invertible. Then there
exists a neighbourhood $U^\prime\times V^\prime$ of $(0,A)$ and a
smooth function $\phi:U^\prime \times V^\prime\to\CN$ with
$\phi(0,A) = 0$, such that if $F(Z,\bar Z, W) = 0$ for some
$(Z,W)\in U^\prime\times V^\prime$, then $Z=\phi (Z,\bar Z,W)$.
Furthermore, for every multiindex $\alpha$, and each $j$, $1\leq
j\leq N$,
\begin{equation}\label{E:impvan}
  D^\alpha \vardop{\phi_j}{Z_k} (Z,\bar Z, W) = 0, \quad 1\leq
  k\leq N,
\end{equation}
if $Z=\phi(Z,\bar Z, W)$, and $\phi$ is  holomorphic in $W$.
Here, $D^\alpha$ denotes the derivative in all the real variables.
\end{thm}
\begin{proof} Let us write $F(Z,\bar Z , W)= F(x,y,W)$ where
$(x,y) \in \RN\times\RN$ are the underlying real coordinates in
$\CN$, as usual identified by $Z_j = x_j + i y_j$. Let us also
choose a neighbourhood $U_0\subset\RN$ of $0$ with the property
that $U_0\times U_0\subset U$. We extend $F$ in the first $2N$
variables almost holomorphically; that is, we have a function
$\tilde F : U_0\times\RN\times U_0\times\RN\times \C^p\to \CN$
with the property that \begin{equation}\label{E:tildef1} \tilde
F(x,x^\prime,y,y^\prime,W)|_{x^\prime = y^\prime = 0} = F(x,y,W)
\end{equation}
and, if we introduce complex coordinates $\xi_k = x_k + i
x_k^\prime$, $\eta_k = y_k + i y_k^\prime$, $1\leq k\leq N$, then
\begin{equation}\label{E:tildef2}
  \left. D^\alpha \vardop{\tilde F_j}{\bar \xi_k} \right|_{x^\prime =
  y^\prime = 0} = \left. D^\alpha \vardop{\tilde F_j}{\bar \eta_k} \right|_{x^\prime =
  y^\prime = 0} = 0, \quad 1\leq j,k\leq N.
\end{equation}
Also, $\tilde F$ is still polynomial in $W$. We introduce new
coordinates $\chi=(\chi_1,\dots, \chi_N) \in\CN$ by
\[\xi_k = \frac{z_k+ \chi_k}{2} , \quad \eta_k = \frac{z_k - \chi_k}{2i}, \quad 1\leq k\leq N, \]
and write $G(Z,\bar Z,\chi,\bar\chi, W) = F(\xi, \bar\xi,
\eta,\bar\eta,W)$. $G$ is smooth in the first $2N$ complex
variables in some neighbourhood of the origin, and polynomial in
$W$. We will now compute the real Jacobian of $G$ with respect to
$Z$ at $(O,A)$. At $(0,A)$, $\vardop{G}{Z}(0,A) =
\vardop{F}{Z}(0,A)$ and $\vardop{G}{\bar Z}(0,A) = 0$, so that we
have
\[ \det \begin{pmatrix}
  \vardop{G}{Z} & \vardop{G}{\bar Z} \\
  \vardop{\bar G}{Z} & \vardop{\bar G}{\bar Z}
\end{pmatrix}(0,A) = \left| \det \vardop{F}{Z} (0,A) \right|^2
\neq 0 \] by assumption. Hence, by the implicit function theorem,
there exists a smooth function $\psi$ defined in some
neighbourhood of $(0,A)$, valued in $\CN$, such that
$Z=\psi(\chi,\bar \chi, W)$ solves the equation $G(Z,\bar Z,
\chi, \bar \chi , W)=0$ uniquely. Here we have already taken into
account that $\psi$ depends holomorphically on $W$, a fact that
the reader will easily check. Since $G(Z,\bar Z, \bar Z, Z
,W)=F(Z,\bar Z,W)$, this implies that if $F(Z,\bar Z,W) = 0$,
then $Z = \psi(\bar Z, Z, W)$.

We let $\phi (Z, \bar Z, W) = \psi (\bar Z, Z, W)$ and claim that
$\phi$ satisfies \eqref{E:impvan}. In fact, computation shows
that $\phi_Z (Z,\bar Z, W)  = \psi_{\bar \chi} = - (G_Z - G_{\bar
Z} \bar G_{\bar Z}^{-1} \bar G_Z )^{-1} (G_{\bar\chi} + G_{\bar
Z} \bar G_{\bar Z}^{-1} \bar G_{\bar \chi})$, where the right
hand side is evaluated at $(\psi(\bar Z, Z, W), \bar \psi(\bar Z,
Z, W), \bar Z, Z, W)$. This formula shows that each $\phi_{j,Z_k}$
is a sum of products each of which contains a factor which is a
derivative of $G$ with respect to $\bar Z$ or $\bar \chi$.

By the definition of $G$, we have that
\[ \vardop{G}{\bar Z} = \frac{1}{2} \vardop{\tilde F}{\bar \xi} +
\frac{1}{2i} \vardop{\tilde F}{\bar \eta}, \quad \vardop{G}{\bar
\chi} = \frac{1}{2} \vardop{\tilde F}{\bar \xi} - \frac{1}{2i}
\vardop{\tilde F}{\bar \eta}.\] By \eqref{E:tildef2} every
derivative of those vanishes if $x^\prime = y^\prime = 0$, which
is in turn the case if $\imag \frac{\phi(\bar Z,Z) + \bar Z}{2} =
0 $ and $\imag \frac{\phi(\bar Z, Z) - Z}{2i} = 0$. But this is
clearly fulfilled if $Z = \phi(\bar Z, Z)$.

The proof is now finished by applying the Leibniz rule, the chain
rule and the observations made above.
\end{proof}

Note that it is clear from the usual implicit function theorem
that we can solve for $N$ of the real variables $(x,y)$. What
this theorem asserts is that we can do so in a special manner.

\section{Proof of Theorem~\ref{T:main}}\label{S:proof}
Let us start by choosing coordinates. There is a neighbourhood
$U$ of $p_0 = 0$ in $\CN$ and a smooth function $\phi:\Cn\times\Rd
\to \Rd$ defined in a neighbourhood $V$ of $0$ such that $M\cap U
= \{ (z,s + i \phi (z,\bar z, s)) \colon (z,s)\in V\}$ with the
property that $\nabla \phi (0) = 0$. Since the conclusion of the
theorem is local, we shall replace $M$ by $M\cap U$, and use this
representation. For suitably chosen open sets $U\subset\Cn$ and
$V\subset\Rd$, consider the diffeomorphism $\Psi: U\times V \to
M$, $\Psi (z,\bar z , s) = (z, s+i \phi( z,\bar z, s))$. We
extend this diffeomorphism almost holomorphically to a map, again
denoted by $\Psi$, from $U\times V\times \Rd$ to $\CN$. $\Psi$ is
a diffeomorphism in an open neighbourhood of $U\times V \times
\{0\}$, and it has the property that for every component $\Psi_l$
of $\Psi$,
\begin{equation}\label{E:almholpsi}
  D^\alpha_{x,y,s,t} \dbar_j \Psi_l (z,s,0) = 0, \quad (z,s)\in U\times V,
\end{equation}
where the derivative is in all the real variables. Equivalently,
\begin{equation}\label{E:almholpsi2}
  D^\alpha_{x,y}D^\beta_s \dbar_j \Psi_l (z,s,0) = O(|t|^\infty), \quad (z,s)\in U\times V,
\end{equation}
uniformly on compact subsets of $U\times V$. That is, for each
$\alpha$, $\beta$, $K\subset U\times V$ compact and every $l\in\N$
there exists a constant $C_l= C_l (\alpha,\beta,K)$ such that
\begin{equation}\label{E:almholpsi3}
  |D^\alpha_{x,y}D^\beta_s \dbar_j \Psi_l (z,s,t)| \leq C_l |t|^l, \quad (z,s)\in
  K.
\end{equation}

We assume that each component $H_j$ of $H$ extends continuously
(and, consequently by a theorem of Rosay \cite{Ro1} already
alluded to above, in a $C^k$-fashion) to a holomorphic function
into a wedge with edge $M$. Let us recall that this means that
with an open convex cone $\Gamma$ in $\Rd$ each $H_j$ extends
continuously to the set $W_\Gamma = \{ Z\in U_0 \colon
\rho(Z,\bar Z) \in \Gamma \}$, where $U_0$ is an open
neighbourhood of $0$ in $\CN$. By choosing $\Gamma$ accordingly,
and possibly shrinking $U_0$, we can in addition assume that each
$H_j$ is continuous and bounded on the closure of $W_\Gamma$, and
in fact smooth up to $b W_\Gamma \setminus M$.

There exists another open, convex cone $\Gamma^\prime$,
relatively closed in $\Gamma$, neighbourhoods $U^\prime\subset U$
and $V^\prime \subset V$ of $0\in\Cn$ and $0\in\Rd$,
respectively, and $\da = (\da_1,\dots,\da_d) >0 $ such that the
wedge $\hat W_{\Gamma^\prime} = \{ (z,s,t) \in U^\prime\times
V^\prime \times \Gamma^\prime \colon 0<t<\da \}$ with flat edge
$U^\prime\times V^\prime$ satisfies $\tilde W_{\Gamma^\prime} =
\Psi (\hat W_{\Gamma^\prime} ) \subset W_\Gamma$. Hence, $h_j =
H_j \circ \Psi$ is well defined on $\hat W_{\Gamma^\prime}$ for
$1\leq j \leq d$, extends continuously to $\hat
W_{\Gamma^\prime}$ and is smooth up to $b \hat W_{\Gamma^\prime}
\setminus U^\prime \times V^\prime$. Since the conclusion of the
theorem is local, we can replace $U$ by $U^\prime$ and $V$ by
$V^\prime$. Furthermore, by shrinking the the neighbourhoods once
more if necessary, we have that there exist positive constants
$C_1$ and $C_2$ such that (here, $d(A,B) $ denotes the distance
between a compact set $A$ and a closed set $B$)
\begin{equation}\label{E:distest1} C_1 d((z,s,t),b\hat W_{\Gamma^\prime}) \leq
d(\Psi (z,s,t), b\tilde W_{\Gamma^\prime} ) \leq C_2
d((z,s,t),b\hat W_{\Gamma^\prime}). \end{equation}

Our next claim is that we can replace $\Gamma^\prime$ by the
standard cone $\Rd_+= \{ t\in\Rd \colon t>0 \}$. In fact, since
$\Gamma^\prime$ is open, we can find $d$ linearly independent
vectors $v_1,\dots, v_j$ in $\Gamma^\prime$. The linear mapping
$T$ which maps $v_j$ to the $j$-th standard basis vector $e_j$ is
invertible, and $T^{-1} (\Rd_+) \subset\Gamma^\prime$ by
convexity. Then we can make a complex linear change of
coordinates by setting $(z^\prime,s^\prime, t^\prime)=(z,T^{-1}s,
T^{-1} t)$. Since this coordinate change is linear and there
exist positive constants $C_1$ and $C_2$ with $C_1 |t| \leq
|t^\prime| \leq C_2 |t|$, \eqref{E:almholpsi},
\eqref{E:almholpsi2}, and \eqref{E:almholpsi3} also hold in the
new coordinates. We need just one more coordinate change.

\begin{cl}\label{claim1} There exists a $\da > 0$, coordinates $(z,s,t)$ and positive
constants $C_1$ and $C_2$ such that $\Psi(z,s,t) \subset \tilde
W_{\Gamma^\prime}$ for $(z,s)\in U\times V$, $0< t <\da$ and $C_1
|t| \leq d(b \tilde W_{\Gamma^\prime}, \Psi(z,s,t) ) \leq C_2
|t|$ for $(z,s)\in U\times V$, $0\leq t\leq \da$.
\end{cl}

\begin{proof} Let $e_j$ denote the $j$-th standard basis vector in $\Rd$,
$1\leq j \leq d$. If $t = (t_1,\dots,t_d) \in \Rd_+$, then
clearly $d(t,b\Rd_+) = \min_{j=1}^d t_j$. For $\eps
> 0$ consider the vectors $v_j = e_j + \eps \sum_{l\neq j} e_l$,
$1\leq j \leq d$. For $\eps$ small enough, these are linearly
independent. We now consider the linear change of coordinates
given by $z^\prime = z$, $t^\prime =
(t_1^\prime,\dots,t_d^\prime) \mapsto \sum_{j=1}^d t_j^\prime
v_j$, $s^\prime = (s_1^\prime,\dots,s_d^\prime) \mapsto
\sum_{j=1}^d s_j^\prime v_j$. By \eqref{E:distest1} it is enough
to show that there exist positive constants $C_1$ and $C_2$ such
that $C_1 |t^\prime| \leq d(t,b\Rd_+) \leq C_2 |t^\prime|$. The
existence of $C_2$ is clear. But if $\eps <1$, then  $d(t,b\Rd_+)
= \min_{j=1}^d t_j = \min_{j=1}^d (t_j^\prime + \eps\sum_{l\neq j}
t_l^\prime) \geq \eps(t_1^\prime+ \dots t_d^\prime)\geq
\frac{\eps}{d} |t|$. An appropriate choice for $\da$ finishes the
argument. \end{proof}

We are going to use the notation introduced in section
\ref{S:boundaryvalues}; that is, we let $\Om_+ = U\times V \times
\{ t\in\Rd_+ \colon 0< t <\da\}$. We let $h_j = H_j\circ \Psi$ on
$\Om_+$.

\begin{cl}\label{claim2} $h_j \in \spaceainf$ for $1\leq j \leq N^\prime$.
\end{cl}

\begin{proof} By all the choices above, $h_j$ satisfies the smoothness
assumptions. Let us first check that every derivative of $h_j$ is
of slow growth. Since $H_j$ is holomorphic in $\tilde
W_{\Gamma^\prime}$ and continuous on its closure, the Cauchy
estimates imply that we have an estimate of the form
\begin{equation}\label{E:cauchy1}
  |\partial^\beta H_j (Z) | \leq C_{\beta} (d(Z,b\tilde
  W_{\Gamma^\prime}))^{-|\beta|}
\end{equation}
for each  $\beta$, where $\partial^\beta$ denotes
$\vardop{^{|\beta|}}{Z^\beta}$. By the chain rule,
$D^\alpha_{x,y,s} h_j (z,s,t)$ is a sum of products of derivatives
of $\Psi$ (which are bounded) and a derivative of $H_j$ with
respect to $Z$, evaluated at $\Psi (z,s,t)$, of order at most
$|\alpha|$. Hence, by \eqref{E:cauchy1} and claim ~\ref{claim1} we
conclude that there exists a positive constant $C$ such that
\begin{equation}\label{E:cauchy2}
  |D^\alpha_{x,y,s} h_j (z,s,t) | \leq C_{\alpha} |t|^{-|\alpha|}.
\end{equation}
We now have to estimate the derivatives of $\dbar_m h_j$ for
$1\leq m\leq d$. But $\dbar_m h_j = \sum_{l=1}^{N^\prime}
\vardop{H_j}{Z_l} \dbar_m \Psi_l$. Hence, if we take an arbitrary
derivative of $\dbar_m h_j$, we get a sum of products of
derivatives of components of $\Psi$ and a derivative of $H_j$ with
respect to $Z$ each of which contains a term of the form $\dbar_m
\Psi_l$. By \eqref{E:cauchy1} and \eqref{E:almholpsi3} we
conclude that for each compact set $K\subset U\times V$ and each
$k\in \N$ there exists a positive constant $C_k$ with
$|D^\alpha_{x,y,s} \dbar_m h_j (z,s,t)| \leq C_k |t|^k$. This
proves claim~\ref{claim2}.
\end{proof}
We now equip $U\times V$ with the CR-structure of $M$; that is, a
basis of the CR-vector fields near $0$ is given by $\Lambda_j =
\Psi^* L_j$ for $1\leq j\leq n$. We almost holomorphically extend
the coefficients of the $\Lambda_j$ to get smooth vector fields
on an open subset of $\Cn\times\Rd\times\Rd$ containing $0$.
\begin{cl}\label{claim3} For each $j$, $1\leq j \leq N^\prime$, there exists a smooth function
$\phi_j (Z^\prime, \bar Z^\prime, W)$ defined in an open
neighbourhood of $(0,(\Lambda^\alpha h (0))_{|\alpha | \leq
k_0})$ in $\CN\times \C^{K(k_0)}$ ($K(k_0)$ denoting $N^\prime
|\{ \alpha \colon |\alpha|\leq k_0 \}|$) such that
\begin{equation}\label{E:reflid}
   h_j (z,s,0) = \phi_j ( h(z,s,0), \overline{h(z,s,0)} , (\overline{\Lambda^\alpha
   h(z,s,0)})_{|\alpha | \leq k_0});
\end{equation}
here, we write $h= (h_1, \dots, h_{N^\prime})$. Furthermore, after
possibly shrinking $U$ and $V$, the right hand side of
\eqref{E:reflid} defines a function in $\spaceam$.
\end{cl}
This last claim of course establishes Theorem~\ref{T:main}; since
$h_j \in \spacea$ by Claim~\ref{claim2} and by Claim~\ref{claim3}
$h_j\in \spaceam$, we can apply Theorem~\ref{T:edgeofthewedge} to
see that $h_j$ is smooth.
\begin{proof}
By the chain rule, we have smooth functions $\Phi_{l,
\alpha}(Z^\prime, \bar Z^\prime, W)$ for $|\alpha| \leq k_0$,
$1\leq l\leq d^\prime$, defined in a neighbourhood of $\{
0\}\times \C^{K(k_0)}$ in $\CN\times\C^{K(k_0)}$, polynomial in
the last $K(k_0)$ variables, such that
\begin{equation}\label{E:thephialpha}
  \Lambda^\alpha \rho^\prime_l (h ,\bar h ) (z,s,0) =
  \Phi_{l,\alpha}
  (h(z,s,0),\overline{h(z,s,0)},(\Lambda^\alpha \overline{
  h (z,s,0)} )_{|\alpha| \leq k_0} ),
\end{equation}
and  $\Lambda^\alpha \rho^\prime_{l,Z^\prime} (h,\bar h) |_0 =
\Phi_{l,\alpha, Z^\prime} (0,0,(\Lambda^\alpha
  h (0,0,0) )_{|\alpha| \leq k_0} ))$. By Definition~\ref{D:nondeg} we
can choose $\alpha^1,\dots, \alpha^{N^\prime}$ and $l^1,\dots,
l^{N^\prime}$ such that if we set $\Phi = (\Phi_{l^1,\alpha^1},
\dots, \Phi_{l^{N^\prime},\alpha^{N^\prime}} )$, then
$\Phi_{Z^\prime} (0)$ is invertible. Hence, we can apply
Theorem~\ref{T:ift}; let us call the solution $\phi$. Then
$\phi_j$ satisfies \eqref{E:reflid}, and we shrink $U$ and $V$
and choose $\da$ in such a way that $g_j(z,s,t) = \phi_j
(h(z,s,-t), \overline{h(z,s,-t)}, (\Lambda^\alpha \overline{
h(z,s,-t)})_{|\alpha|\leq k_0})$ is well defined and continuous
in a neighbourhood of $\bar \Om_-$. It is easily checked that
$g_j$ is a function in $\spaceam$ as a consequence of
\eqref{E:impvan} and the fact that each $h_j\in \spaceainf$.
First note that this implies $\overline{h_j (z,s,-t)}
\in\spaceaminf$, and by Lemma~\ref{L:difflemma}, $\Lambda^\alpha
\overline{ h_j (z,s,-t)}\in\spaceaminf$ for each $\alpha$. Now,
each derivative $D^\beta$ of $g_j$ is a sum of products of
derivatives of $\phi_j$ (which are uniformly bounded on $\Om_-$)
and derivatives of $h$, $\bar h$, and $\Lambda^\alpha \bar h$,
all of which fulfill the analog of \eqref{E:slowgrowth} on
$\Om_-$. So $g_j$ fulfills the analog of \eqref{E:slowgrowtha} on
$\Om_-$. Next, we compute the derivative of $g_j$ with respect to
$\bar w_k$. We have that
\[\vardop{g_j}{\bar w_k} = \sum_{l=1}^{N\prime} \vardop{\phi_j}{Z^\prime_l}
\vardop{h_l}{\bar w_k} + \sum_{l=1}^{N\prime} \vardop{\phi_j}{\bar
Z^\prime_l} \vardop{\bar h_l}{\bar w_k} + \sum_{|\alpha| \leq k_0}
\vardop{\phi}{W_\alpha} \vardop{\Lambda^\alpha \bar h}{\bar w_k}.
\]

Applying any derivative $D^\beta$, we see that the first sum
gives rise to products of derivatives of
$\vardop{\phi_j}{Z^\prime_l}$ and derivatives of $h$, $\bar h$,
and $\Lambda^\alpha \bar h$. Now the derivatives of $\phi_j$
fulfill \eqref{E:impvan}. Since on $t = 0$, $h = \phi (h,\bar h,
(\Lambda^\alpha \bar h)_{|\alpha |\leq k_0} )$, we conclude that
$h - \phi (h,\bar h, (\Lambda^\alpha \bar h)_{|\alpha |\leq k_0}
) = O(|t|)$. But by \eqref{E:impvan}, any derivative of
$\vardop{\phi_j}{Z_l} (Z,\bar Z, W)$ is $O(|Z - \phi(Z,\bar Z,
W)|^\infty)$, so that derivatives of
$\vardop{\phi_j}{Z^\prime_l}$ evaluated at $(h,\bar h,
(\Lambda^\alpha \bar h)_{|\alpha |\leq k_0})$ are
$O(|t|^\infty)$. All the other terms in the product are
$O(|t|^{-s})$ for some $s$, so that the terms coming from the
first sum are actually $O(|t|^\infty)$. For the second and third
sum, a similar argument using that $\bar h$ and $\Lambda^\alpha
\bar h$ are in $\spaceaminf$ implies that all the terms arising
from them are $O(|t|^\infty)$. All in all, we conclude that $g_j
\in \spaceaminf$, which finishes the proof.
\end{proof}
\bibliographystyle{plain}
\bibliography{smoothdeg}
\end{document}